\documentclass[runningheads]{llncs}

\usepackage{amssymb}
\usepackage{url}

\begin{document}

\mainmatter  

\title{ An Alternative Method for Primary Decomposition of Zero-dimensional Ideals over Finite Fields }

\author{Yongbin Li}

\institute{School of Mathematical Sciences\\
      University of Electronic Science and Technology of China\\
       Chengdu, Sichuan 611731,  China\\
       Email:yongbinli@uestc.edu.cn}

%
%

\maketitle

\begin{abstract}
 We present an alternative method for computing primary decomposition of
 zero-dimensional ideals over finite fields.
Based upon the further decomposition of the invariant subspace of
the Frobenius map acting on the quotient algebra in the algorithm
given by S. Gao, D. Wan and M. Wang in 2008, we get an alternative
approach to compute all the primary components at once. As one
example of our method, an improvement of Berlekamp's algorithm by
theoretical considerations which computes the factorization of
univariate polynomials over finite fields is also obtained.

\end{abstract}

\section{Introduction}
Let $k$ be a  field   and $k[x_1,\ldots,x_n]$ (or $k[{\bf x}]$ for
short) be the ring of polynomials in the variables
$x_1,\ldots,x_n$ with coefficients in $k$. An ideal $I\subset
k[{\bf x}]$ is called  zero-dimensional if the quotient algebra
$k[{\bf x}]/I$ is a finite-dimensional $k-$vector space.
 There are several well-know algorithms for computing primary
ideal decomposition based on zero-dimensional decomposition, and
we refer the readers to \cite{CM02,NKY,Shi96,GT88}.

Gao, Wan and Wang in \cite{GWM2009} present an interesting
approach to compute primary decomposition of zero-dimensional
ideals over finite fields.
 The method is based on  the invariant subspace of the
Frobenius map acting on the quotient algebra $k[{\bf x}]/I$. Since
the dimension of the invariant subspace just equals the number of
primary components, a basis of the invariant subspace leads a
complete primary decomposition of $I$ by computing  Gr\"obner
bases.

In the method in \cite{GWM2009}, if one chooses  an element of the
basis of the invariant subspace which  is  {\it separable } for
$I$, then  all the primary components can be computed at once.
Otherwise, the further decomposition is necessary even though the
probability of separable element in the invariant subspace is not
low in most cases, see Proposition 3.2 in \cite{GWM2009} for
details. In this paper, we aim to find an approach to  decompose
the above invariant subspace into a direct product of several
one-dimensional $k-$algebras in Proposition 1 in Section 3. Based
upon this theoretical work, we get  an alternative approach to
compute primary decomposition of zero-dimensional ideals over
finite fields, which allows us to find  all the primary components
completely.

\section{Preliminaries}

In this section, we assume $R$ is a commutative ring. For basic
notations in commutative algebra, we refer the readers to the
monographs \cite{Cox98,Martin:00}.

{\bf Definition 1.$^{[5]}$} Given a commutative ring $R$, let $I$
be an ideal in $R$, let ${\bf U}$ be an $R-$module. Then the set
$\{m\in {\bf U}|m\cdot s=0 ,\;\forall\; s\in I\}$ is an
$R-$submodule of ${\bf U}$. It is called the {\it colon module} of
${\bf 0}$ by $I$ in ${\bf U}$, denoted by $M(I)$.

{\bf Lemma 1.} Let $ R$ be a commutative ring with identity
element, and ${\bf U}$ be an  $R-$module. If $I_1,\ldots, I_t$ are
pairwise comaximal ideals in $ R$, then
$$M(I)=M(I_1)\oplus \cdots \oplus M(I_t),$$ where $I=I_1\cap\cdots \cap I_t$ or $I_1\cdots
I_t.$

{\it Proof:}\/ We give a proof by induction on $t$. When $t=2$,
first we prove the following claim
$$M(I)=M(I_1\cap I_2)=M(I_1)+M(I_2).$$
It is obvious that $M(I)=M(I_1\cap I_2)\supseteq M(I_1)+M(I_2) $
by the definition of $M(I)$ and the fact that  $I=I_1\cap
I_2=I_1I_2.$

Since $I_1$ and $I_2$ are pairwise comaximal ideal in $R$, there
exist some $a_1 \in I_1$ and $a_2\in I_2$  such that $a_1+a_2=1$
where $1$ is the identity element in $R,$  which implies that
$m=a_1\cdot m+a_2\cdot m$ for any $m \in M(I_1 \cap I_2)$. The
other direction follows from the fact that $a_1\cdot m \in M(I_2)$
and $a_2 \cdot m\in M(I_1)$.

Since $a_1\cdot m=a_2\cdot m=0$ for any $m\in M(I_1)\cap M(I_2)$,
which means $m=a_1\cdot m+a_2\cdot m=0$. Thus,
$$M(I)=M(I_1)\oplus  M(I_2).$$

Now let $t>2$, consider that $I_1\cap \cdots \cap I_{t-1}$ and
$I_t$ are pairwise comaximail ideals in $R$. It follows from the
induction hypothesis that
$$M(I_1\cap \cdots \cap I_{t-1}\cap I_t)=M(I_1\cap \cdots \cap I_{t-1})\oplus M(
I_t)=M(I_1)\oplus\cdots \oplus M(I_{t-1})\oplus M( I_t).$$ This
completes the proof. $\square$

 We need the following ring-theoretic version
of the Chinese Remainder Theorem in our discussion.

{\bf Lemma 2(Chinese Remainder Theorem).$^{[5]}$} Let $ R$ be a
commutative ring with identity element. If $I_1,\ldots, I_t$ are
pairwise comaximal ideals in $ R$, then the canonical map is  an
isomorphism of R-modules
$$R/I\cong R/I_1\oplus \cdots \oplus R/I_t,$$ where $I=I_1\cap\cdots \cap
I_t$.

Let $k$ be any field, now we consider the above ring
$R=k[x_1,\dots,x_n]$ or $k[{\bf x}]$.

{\bf Theorem 1.} Assume that $I_1,\ldots, I_t$ are pairwise
comaximal ideals in $k[x_1,\ldots,x_n]$, and let $I=I_1\cap\cdots
\cap I_t$ or $I_1\cdots I_t$. Then the canonical map $\Phi$ is an
isomorphism of $k[{\bf x}]-$modules, i.e.,
$$k[{\bf x}]/I=J_1/I\oplus\cdots \oplus J_t/I\cong k[{\bf x}]/I_1\oplus \cdots \oplus k[{\bf x}]/I_t$$
where $J_i=I_1\cap \cdots \cap I_{i-1}\cap I_{i+1}\cap \cdots \cap
I_t$ for $i=1,2,\cdots,t.$ Moreover,  the restriction of $\Phi$ to
$J_i/I$ is an isomorphism of $k[{\bf x}]-$modules
$$J_i/I \cong k[{\bf x}]/I_i,$$
for each $i=1,2,\dots,t.$

{\it Proof:}\/ Consider the $k[{\bf x}]-$ module ${\bf U}=k[{\bf
x}]/I.$ By Lemma 1, we have $$M(I)=M(I_1)\oplus\cdots \oplus
M(I_t).$$
  Since $I_i$ and $J_i$ are comaximal for
each $i$ and
  $M(I)=k[{\bf
x}]/I$ and $M(I_i)=J_i/I$,  it is easy to see that $k[{\bf
x}]/I=J_1/I\oplus \cdots \oplus J_t/I$.

Applying Lemma 2, the  canonical map $\Phi$ is  an isomorphism of
$k[{\bf x}]-$modules
$$k[{\bf x}]/I\cong k[{\bf x}]/I_1\oplus \cdots \oplus k[{\bf x}]/I_t.$$

Furthermore, it implies that each restriction of $\Phi$ to $J_i/I$
is an isomorphism $J_i/I \cong k[{\bf x}]/I_i$ by the proof of
Chinese Remainder Theorem. For the details, please see  the proof
of Lemma 3.7.4 in \cite{Martin:00}. This completes the proof.
$\square$

Let $k$ be any field containing a finite field ${\mathbb F}_q$ as
a subfield. An ideal $I\subset k[{\bf x}]$ is called {\it primary}
if each no-zero zerodivisor of  $k[{\bf x}]/ I$ is a non-zero
nilpotent element. Further $I$ is called  {\it quasi-primary} if
$\sqrt{I}$ is a prime ideal, that is, if $I$ has only one minimal
component and all other components are embedded.

{\bf Definition 2.} Let $k$ be any field containing ${\mathbb
F}_q$ as a subfield, and $I$ be  an ideal in $k[{\bf x}]$, we
establish
  the following ${\mathbb F}_q-$linear transformation
$\Psi_{I}$ from ${\mathbb F}_q-$vector space $k[{\bf x}]/I,$
 to itself, defined by
$$ \Psi_I (\bar{f})= {\bar f}^q-{\bar f}$$ for each ${\bar f}\in k[{\bf
x}].$

In fact $\mathrm{Ker}(\Psi_{I})$ is the invariant subspace of the
Frobenius map acting on $k[{\bf x}]/I$ which plays an essential
role
 in \cite{GWM2009} and our improvement. With the
notation in Definition 2,  Lemma 2.1 in \cite{GWM2009} can be
described as follows.

{\bf Lemma 3.$^{[6]}$} Let $k$ be any field containing ${\mathbb
F}_q$ as a subfield, and $I_0\subset k[{\bf x}]$ be a
quasi-primary ideal. Then
$$\mathrm{Ker}(\Psi_{I_0})={\mathbb
F}_q.$$

Now consider an arbitrary ideal $I\subseteq k[{\bf x}]$. Suppose
$I$ has an irredundant primary decomposition
\begin{eqnarray}
\label{2}
 I=I_1\cap I_2 \cap \cdots \cap I_t,
\end{eqnarray}
where $I_i\in k[{\bf x}]$ are primary ideals, and  $I_i$'s are
pairwise comaximal.

The following result is another version of  Theorem 2.2 in
\cite{GWM2009}. Here we present an alternative proof.

 {\bf
Theorem 2.} Let  $I\subset {\mathbb F}_q[{\bf x}]$ be a
zero-dimensional ideal with $t$ irredundant primary components
$I_1,I_2,\ldots,I_t$. If we consider ${\mathbb F}_q[{\bf x}]/I$ as
${\mathbb F}_q-$vector space,  then
 $\mathrm{Ker}(\Psi_{I})\cong {\mathbb F}_q^t.$

{\it Proof:}\/ From Theorem 1, we know that there exists an
isomorphism of ${\mathbb F}_q-$vector spaces
$${\mathbb F}_q[{\bf x}]/I=J_1/I\oplus \cdots \oplus J_t/I\cong {\mathbb F}_q[{\bf x}]/I_1\oplus \cdots \oplus {\mathbb F}_q[{\bf x}]/I_t$$
and $J_i/I\cong {\mathbb F}_q[{\bf x}]/I_i$ for each $i.$ By Lemma
3, we have $\mathrm{Ker}(\Psi_{I})\cong {\mathbb F}_q^t.$ This
completes the proof. $\square$

\section{Main Results}

Let $I$ be a zero-dimensional ideal in ${\mathbb F}_q[{\bf x}]$.
In the following, we assume that a Gr\"obner basis for $I$ is
already known or computed for certain term order. Then one can
easily find a linear basis for ${\mathbb F}_q[{\bf x}]/ I$ over
${\mathbb F}_q$ by Macaulay¡¯s Basis Theorem in
\cite{Cox98,Martin:00}.

Based upon the result of Theorem 2.2 in \cite{GWM2009} or  Theorem
2 in this paper, we are ready  to decompose
$\mathrm{Ker}(\Psi_{I})$ into a direct product of some
one-dimensional algebras using the following method. Furthermore,
by applying Gr\"obner basis theory, an improving approach to
compute all the primary components of $I$  is given.

 {\bf Lemma 4.}  In the
situation of Theorem 2,  $\mathrm{Ker}(\Psi_{I})$ is a subring of
${\mathbb F}_q[{\bf x}]/I$ and has no non-zero nilpotent elements.

{\it Proof:}\/ It is easy to check that $\mathrm{Ker}(\Psi_{I})$
is a subring of ${\mathbb F}_q[{\bf x}]/I$. We proceed to show
that $\mathrm{Ker}(\Psi_{I})$ has no non-zero nilpotent elements.

Suppose there is some $\bar{f}\in \mathrm{Ker}(\Psi_{I})$ and a
positive integer $m$ satisfying $\bar{f}^m=\bar{0}.$ Consider that
the greatest common divisor of  polynomials $x^m$ and $ x^q-x $ is
$x$ in $ {\mathbb F}_q[x] $, there exist some $u(x),v(x)\in
{\mathbb F}_q[x] $ such that $x=u(x)x^m+v(x)(x^q-x)$. It implies
that
$$\bar{f}=u(\bar{f})\bar{f}^m+v(\bar{f})(\bar{f}^q-\bar{f})=\bar{0}$$
by $\bar{f}^q=\bar{f}$. This completes the proof. $\square$

 The proof of the following theorem presents an approach to
 decompose a finite
dimensional $k_0-$algebra which has no non-zero nilpotent elements
.

 {\bf Proposition 1.} Let $k_0$ be a field, and ${\bf V}$ be a
 $k_0-$algebra with $\mathrm{dim}_{k_0}({\bf V})=t$. If
${\bf V}$ has no non-zero nilpotent elements, then there exist
$\bar{g}_1,\ldots,\bar{g}_t\in {\bf V}$ such that ${\bf
V}=\mathrm{span}_{k_0}\langle \bar{g}_1,\ldots,\bar{g}_t\rangle $
and $\bar{g}_i\bar{g}_j=\bar{0}$ for any $i\ne j.$ Furthermore,
${\bf V}$ can be written as  a direct product of several
one-dimensional $k_0-$algebras
$$ {\bf V}=\langle \bar{g}_1\rangle \oplus \cdots \oplus \langle\bar {g}_t \rangle.$$

{\it Proof:}\/ Let $\{\bar{f_1},\cdots,\bar{f_t}\}$ be a basis of
the $k_0-$vector space ${\bf V}$. For a given no-zero  zerodivisor
$\bar{h}\in {\bf V}$, we define that  next ${\mathbb F}_q-$linear
transformation
\begin{eqnarray*}
{\bf V} & \longrightarrow & {\bf V},\\
 \phi_{\bar{h}}:\; \bar{g} & \longrightarrow & \bar{h}\bar{g}.
\end{eqnarray*}
It is well known that
$\mathrm{dim}_{k_0}(\mathrm{Ker}(\phi_{\bar{h}}))+\mathrm{dim}_{k_0}(\mathrm{Im}(\phi_{\bar{h}}))=t.$

We claim that $${\bf V}=\mathrm{Ker}(\phi_{\bar{h}})\oplus
\mathrm{Im}(\phi_{\bar{h}}).$$ To prove the claim, it suffices to
show $\mathrm{Ker}(\phi_{\bar{h}})\cap
\mathrm{Im}(\phi_{\bar{h}})=\{\bar{0}\}.$ For any $\bar{p}\in
\mathrm{Ker}(\phi_{\bar{h}})\cap \mathrm{Im}(\phi_{\bar{h}})$, it
implies that there exists some $\bar{q}\in {\bf V}$ such that
$\bar{p}=\bar{h}\bar{q}$. Notice that
$\bar{h}\bar{p}=\bar{h}^2\bar{q}=\bar{0}$, hence
$\bar{p}^2=\bar{0}.$ Since ${\bf V}$ has no non-zero nilpotent
elements, we have  $\bar{p}=\bar{0}$. This  finishes the proof of
our claim. In addition both $\mathrm{Ker}(\phi_{\bar{h}})$ and
$\mathrm{Im}(\phi_{\bar{h}})$ are $k_0-$algebras. Thus the claim
implies a direct product of $k_0-$algebras
$${\bf V}=\mathrm{Ker}(\phi_{\bar{h}})\oplus
\mathrm{Im}(\phi_{\bar{h}}).$$

With an analogous operation $\mathrm{Ker}(\phi_{\bar{h}})$ and $
\mathrm{Im}(\phi_{\bar{h}})$, respectively, we can decompose ${\bf
V}$ into a direct product of one-dimensional $k_0-$algebras
$$ {\bf V}=\langle\bar{g}_1\rangle\oplus \cdots \oplus \langle\bar {g}_t\rangle.$$
This completes the proof. $\square$

{\bf Remark:}  Since the number of regular elements of ${\bf V}$
is $|k_0|$, most of elements are non-zero divisors.

Based upon the above results, we can decompose
$\mathrm{Ker}(\Psi_{I})$ in Theorem 2.

{\bf Theorem 3.} Let ${\mathbb F}_q$ be a finite field. Suppose
$I\subset {\mathbb F}_q[{\bf x}]$ is a zero-dimensional ideal with
$t$ irredundant primary components $I_1,I_2,\ldots,I_t$. Then
there exists a  direct product of one dimensional ${\mathbb
F}_q-$algebras
$$ \mathrm{Ker}(\Psi_{I})=\langle\bar{h}_1\rangle\oplus \cdots \oplus \langle\bar {h}_t\rangle,$$
with each $h_i\in J_i \setminus I_i$ and $\bar{h}_i^2=\bar{h}_i$
where $J_i=I_1\cap \cdots \cap I_{i-1}\cap I_{i+1}\cap \cdots \cap
I_t$ for $i=1,2,\dots,t.$

{\it Proof:}\/ We regard  $ {\mathbb F}_q[{\bf x}]/I$ as an
${\mathbb F}_q-$algebra. With the notation in Definition 2, let
${\bf V}_0=\mathrm{Ker}(\Psi_{I})$. It follows from Theorem 2 that
${\bf V}_0$ is a finite dimensional ${\mathbb F}_q-$algebra with
$\mathrm{dim}_{{\mathbb F}_q}({\bf V}_0)=t$. Furthermore, it
follows that ${\bf V}_0$ has no non-zero nilpotent elements from
Lemma 4.

Applying Proposition 1, one can get the next direct product of one
dimensional ${\mathbb F}_q-$algebras
$$ {\bf V}_0=\langle\bar{g}_1\rangle\oplus \cdots \oplus \langle\bar {g}_t\rangle.$$
Taking some permutation of $\bar{g}_1,\dots,\bar{g}_t$, we can
make sure that $g_i\in J_i$ for $i=1,2, \dots, t.$ Since there
exists some $k_i\in {\mathbb F}_q \setminus \{0\} $ satisfying
$\bar{g}_i^2=k_i\bar{g}_i$  for each $i$, we can easily get some
$\bar{h}_i\in \langle \bar{g}_i\rangle\setminus \{0\} $ such that
$\bar{h}_i^2=\bar{h}_i$ for $i=1,2,\dots, t.$
 This completes
the proof. $\square$

The following example which is the same one given in
\cite{GWM2009} illustrates our decomposition.

{\bf Example 1.} Consider the ideal

$$I=\langle y^2-xz,z^2-x^2y,x+y+z-1 \rangle \subset {\mathbb F}_5[x,y,z].$$
Under the lex order with $x\succ y\succ z$, $I$ has a Gr\"obner
basis
\begin{eqnarray*}
G& =& [x+y+z-1,y^2+3y-2z^4+z^3+2z^2+z,\\
 & & yz+2y+2z^4-z^3-z^2-2z,z^5-z^4+3z^3-z^2+2z],
\end{eqnarray*}
Let $R={\mathbb F}_5[x,y,z]/I$. By Macaulay¡¯s Basis Theorem, we
can get a ${\mathbb F}_5-$basis $B$ of $R$ with
$B=\{z^4,z^3,z^2,z,y,1\}.$ Using  the matrix of ${\mathbb
F}_5-$linear transformation $\Psi_{I}$ on the basis $B$, or
referring to \cite{GWM2009}, one can  easily compute that
$${\bf V}=\mathrm{Ker}({\Psi_{I}})=\{\bar{g}\in R\;|\; \bar{g}^5=\bar{g}\}=\mathrm{span}_{{\mathbb F}_5}\langle \bar{g}_1,\bar{g}_2,\bar{g}_3,\bar{g}_4 \rangle$$
where
$$\bar{g}_1=1,\;\bar{g}_2=z-z^2,\;\bar{g}_3=z^2+z^3,\;\bar{g}_4=z^3-2z^4.$$
We proceed to show how to decompose ${\bf V}$ using the method
given in the proof of  Proposition 1.

As $\bar{h}=\bar{g}_4\in {\bf V}$ is a zerodivisor, we define the
next ${\mathbb F}_5-$linear transformation
\begin{eqnarray*}
{\bf V} & \longrightarrow & {\bf V},\\
 \phi_{\bar{g}_4}:\; \bar{g} & \longrightarrow & \bar{g}_4\bar{g}.
\end{eqnarray*}
The  matrix of $\phi_{\bar{g}_4}$ under the basis
$\bar{g}_1,\bar{g}_2,\bar{g}_3,\bar{g}_4$ of ${\bf V}$ is

$$A_{\bar{g}_4}=\left[ \begin {array}{cccc}
0&0&0&0\\\noalign{\medskip}0&2&2&-1\\\noalign{\medskip}0&2&2&-1\\\noalign{\medskip}1&3&2&1\end
{array}
 \right]$$ such that
 $(\phi_{\bar{g}_4}(\bar{g}_1),\phi_{\bar{g}_4}(\bar{g}_2),\phi_{\bar{g}_4}(\bar{g}_3),\phi_{\bar{g}_4}(\bar{g}_4))
 =(\bar{g}_1,\bar{g}_2,\bar{g}_3,\bar{g}_4)A_{\bar{g}_4}.$ We can
 easily compute that
 $$\mathrm{Ker}(\phi_{\bar{g}_4})=\mathrm{span}_{{\mathbb
 F}_5}\langle \bar{h}_1,\bar{h}_2\rangle\;\mathrm{and}\; \mathrm{Im}(\phi_{\bar{g}_4})=\mathrm{span}_{{\mathbb
 F}_5}\langle\bar{h}_3,\bar{h}_4\rangle $$ where
\begin{eqnarray*}
\bar{h}_1& =& -4-\bar{g}_2+\bar{g}_3=-4-z+2z^2+z^3,\\
 \bar{h}_2& =& 3\bar{g}_2+\bar{g}_4=3z-3z^2+z^3-2z^4,\\
 \bar{h}_3& =& \bar{g}_2+\bar{g}_3=z^3+z,\\
 \bar{h}_4& =& \bar{g}_4=z^3-2z^4.\\
\end{eqnarray*}
Let ${\bf V}_1=\mathrm{Ker}(\phi_{\bar{g}_4})$ and ${\bf
V}_2=\mathrm{Im}(\phi_{\bar{g}_4})$.
 It is easy to check that both ${\bf V}_1$ and ${\bf V}_2$ are
 two-dimensional  ${\mathbb F}_5-$algebras and  have no
 non-zero nilpotent elements. Moreover, there is a  direct product of
 two-dimensional ${\mathbb F}_5-$algebras
 $${\bf V}={\bf V}_1\oplus {\bf V}_2.$$

 Continuing this way, we can compute analogously as the following
 direct product of
 one-dimensional ${\mathbb F}_5-$algebras
$${\bf V}_1=\langle\bar{h}_2\rangle\oplus \langle\bar{h}_1+\bar{h}_2\rangle,\;\;{\bf V}_1=\langle\bar{h}_3\rangle\oplus
\langle\bar{h}_3+\bar{h}_4\rangle.$$ Thus, ${\bf
V}=\langle\bar{h}_2\rangle\oplus
\langle\bar{h}_1+\bar{h}_2\rangle\oplus\langle\bar{h}_3\rangle\oplus
\langle\bar{h}_3+\bar{h}_4\rangle$ such that
\begin{eqnarray*}
(\bar{h}_2)^2& =& 2\bar{h}_2,\;\;
 (\bar{h}_1+\bar{h}_2)^2 =\bar{h}_1+\bar{h}_2,\\
 (\bar{h}_3)^2& =& \bar{h}_3,\;\;
 (\bar{h}_3+\bar{h}_4)^2 = \bar{h}_3+\bar{h}_4.
 \end{eqnarray*}

Furthermore, since $ (3\bar{h}_2)^2 = 3\bar{h}_2$, we have $${\bf
V}=\langle3\bar{h}_2\rangle\oplus
\langle\bar{h}_1+\bar{h}_2\rangle\oplus\langle\bar{h}_3\rangle\oplus
\langle\bar{h}_3+\bar{h}_4\rangle.$$

It remains to compute  all the primary components of $I$ by the
following results.

{\bf Proposition 2.} Let $ R$ be a commutative ring with identity
element 1. If there exist $f_0\in I_0,\;g_0\in J_0$ such that
$f_0+g_0=1$, namely, $I_0, J_0$ are two comaximal ideals in $ R$,
 then
$$I:\langle g_0\rangle^{\infty}=I_0,$$ where $I=I_0\cap J_0.$

{\it Proof:}\/ We first show that $I:\langle
g_0\rangle^{\infty}\supseteq I_0.$ Given any $h \in I_0$, we have
$h=hf_0+hg_0.$  It implies that $hg_0=h-hf_0\in I_0$. Since
$hg_0\in J_0$, $hg_0\in I_0\cap J_0=I$. Thus, $h\in I:\langle
g_0\rangle\subseteq I:\langle g_0\rangle^{\infty}.$

In another direction, for any $h\in I:\langle
g_0\rangle^{\infty}$, it means that there exists some integer
$k>0$ such that
$$hg_0^k\in I.$$
According to $(f_0+g_0)^k=1$, we have that there is some $f_0^*\in
I_0$ such that $f_0^*+g_0^k=1$. Hence
$$h=hf_0^*+hg_0^k\in I_0+I \subseteq I_0.$$
This completes the proof. $\square$

One can compute $I:\langle g_0\rangle^{\infty}$ in $k[{\bf x}]$ by
the following Gr\"obner basis theory, see \cite{Cox:06,Martin:00}
for the details.

{\bf Lemma 5.} Let $k$ be a field and $I$  an ideal generated by
$\{f_1,\ldots,f_s\}$ in $k[{\bf x}]$, and some  $g_0\in k[{\bf
x}]$. If $G^*$ is the Gr\"obner basis of
$\{f_1,\ldots,f_s,1-ug_0\}$ in $k[{\bf x},u]$ with respect to the
purely lexicographical order determined with $x_i\prec u$. Then
$$I:\langle g_0\rangle^{\infty}=\langle G^*\cap k[{\bf x}]\rangle.$$

Next we proceed to present an alternative algorithm for computing
primary decomposition of
 zero-dimensional ideals over finite fields .
 \vspace{.1in}

\vspace{.1in} \noindent {\bf  Primary Decomposition:} $I_1,\dots ,
I_t \leftarrow I$. Given a zero-dimensional ideal in ${\mathbb
F}_q[\bf {x}]$, this algorithm computes an irredundant primary
decomposition $I=I_1\cap \cdots \cap I_t$.

\begin{description}
  \item[D1.]Computing a Gr\"obner basis $G$ of $I$,
   get the  basis $B$ for ${\mathbb F}_q[{\bf x}]/I$ by Macaulay¡¯s Basis Theorem.

\item[D2.]Compute $\mathrm{Ker}(\Psi_{I})$ and $t$:

\begin{description}
  \item[D2.1.]Compute $\mathrm{Ker}(\Psi_{I})$ by the matrix of ${\mathbb F}_q-$linear transformation $\Psi_{I}$ on the basis $B$ and
 set
       ${\bf V} \leftarrow \mathrm{Ker}(\Psi_{I})$.
   \item[D2.2.]Set $t \leftarrow \mathrm{dim}_{{\mathbb F}_q}({\bf V})$.

  \end{description}
  \item [D3.]Decompose ${\bf V}$ into a  direct product of one-dimensional ${\mathbb F}_q-$algebras by the method given in the proof
  of
  Proposition 1,
   $$ {\bf V}=\langle \bar{h}_1\rangle\oplus \cdots \oplus \langle \bar {h}_t\rangle.$$

   \item [D4.]Compute each Gr\"obner basis $G^*_i$ of  $I \cup \{1-u\bar{h}_i\}$ in ${\mathbb F}_q[{\bf x},u]$ by
   Lemma 5. Set $I_i \leftarrow \langle G^*_i\cap {\mathbb F}_q[{\bf x}]\rangle$ for
   $i=1,2,\dots,t.$

\end{description}

{\bf Example 2.}Continued from Example 1,
 compute Gr\"obner basis $G^*_1$ of  $\{y^2-xz, z^2-x^2y, x+y+z-1,1-u\bar{h}_1\}$ in ${\mathbb F}_5[x,y,z,u]$ by
   Lemma 5.
   $G_1^*=[z, y, x+4, 4+u],$ so $$ I_1=\langle z, y, x+4\rangle.$$

 Similarly, we can compute that
 \begin{eqnarray*}
 G_2^*&=& [z+2, y^2+3y+1, x+y+2, 2+u],\\
 G_3^*&=& [z^2+4z+2, y+2z+1, x+4z+3, 4+u],\\
 G_4^*&=& [z+3, y+4, x+2, 4+u].
 \end{eqnarray*}
Therefore,
 \begin{eqnarray*}
 I_2&=& \langle z+2, y^2+3y+1, x+y+2\rangle,\\
 I_3&=& \langle z^2+4z+2, y+2z+1, x+4z+3\rangle,\\
 I_4&=& \langle z+3, y+4, x+2\rangle.
 \end{eqnarray*}

From the above , we have the following irredundant primary
decomposition
$$I=I_1\cap I_2\cap I_3 \cap I_4.$$

As a direct application of our approach, we easily give to an
alternative method  of Berlekamp's algorithm which computes the
factorization of univariate polynomials over finite fields, see
\cite{Berlek67} for the details.

{\bf Example 3.} Consider the polynomial

$$f={x}^{6}+{x}^{5}+{x}^{4}+2 \in{\mathbb F}_3[x].$$
 Let $B=\{x^5,x^4,x^3,x^2,x,1\}$ which is  an ${\mathbb F}_3-$basis of ${\mathbb F}_3[x]/\langle f\rangle$.

Then we have the following ${\mathbb F}_3-$linear transformation
$$\Psi_{\langle f\rangle}\;:\; {\mathbb F}_3[x]/\langle f\rangle \rightarrow {\mathbb
F}_3[x]/\langle f\rangle$$  defined by
$$ \Psi_{\langle f\rangle} (\bar{g})= {\bar g}^3-{\bar g}$$ for ${\bar g}\in {\mathbb
F}_3[x]/\langle f\rangle.$

 The matrix of ${\mathbb
F}_3-$linear transformation $\Psi_{\langle f\rangle}$ on the basis
$B$ is obtained as follows:
$$\left[ \begin {array}{cccccc} 0&0&2&2&0&0\\\noalign{\medskip}0&0&2&2&0&0\\\noalign{\medskip}2&1&0&0&1&0\\\noalign{\medskip}0&2&2&2&0&0
\\\noalign{\medskip}1&0&0&0&2&0\\\noalign{\medskip}0&2&1&1&0&0
\end {array} \right]
.$$

We can  easily compute that
$${\bf V}=\mathrm{Ker}({\Psi_{\langle f\rangle}})=\{\bar{g}\in {\mathbb F}_3[x]/\langle f\rangle\;|\; \bar{g}^3=\bar{g}\}=\mathrm{span}_{{\mathbb F}_3}\langle\bar{g_1},\bar{g_2},\bar{g_3}\rangle$$
where
$$\bar{g}_1=1,\;\bar{g}_2=-x^3+x^2,\;\bar{g}_3=x^5+x.$$
We proceed to show how to decompose ${\bf V}$ by the method given
in the proof of  Proposition 1.

As   $\bar{g}_3\in {\bf V}$ is a zerodivisor, we suppose that the
next ${\mathbb F}_3-$linear transformation
\begin{eqnarray*}
{\bf V} & \longrightarrow & {\bf V},\\
 \phi_{\bar{g}_3}:\; \bar{g} & \longrightarrow & \bar{g}_3\bar{g}.
\end{eqnarray*}
The  matrix of $\phi_{\bar{g}_3}$ under the basis
$\{\bar{g}_1,\bar{g}_2,\bar{g}_3\}$ of ${\bf V}$ is

$$A_{\bar{g}_3}=\left[ \begin {array}{ccc} 0&-1&1\\\noalign{\medskip}0&-1&1\\\noalign{\medskip}1&-1&1\end {array} \right]
$$ such that
 $(\phi_{\bar{g}_3}(\bar{g}_1),\phi_{\bar{g}_3}(\bar{g}_2),\phi_{\bar{g}_3}(\bar{g}_3))
 =(\bar{g}_1,\bar{g}_2,\bar{g}_3)A_{\bar{g}_3}.$ We can
readily compute that
 $$\mathrm{Ker}(\phi_{\bar{g}_3})=\mathrm{span}_{{\mathbb
 F}_3}\langle\bar{s}_1\rangle\;\mathrm{and}\; \mathrm{Im}(\phi_{\bar{g}_3})=\mathrm{span}_{{\mathbb
 F}_3}\langle\bar{s}_2,\bar{s}_3\rangle$$ where
\begin{eqnarray*}
\bar{s}_1& =& \bar{g}_2+\bar{g}_3=x^5+2x^3+x^2+x,\\
 \bar{s}_2& =& \bar{g}_3=x^5+x,\\
 \bar{s}_3& =& \bar{g}_1+ \bar{g}_2+\bar{g}_3=x^5+2x^3+x^2+x+1.\\
 \end{eqnarray*}
Let ${\bf V}_1=\mathrm{Ker}(\phi_{\bar{g}_3})$ and ${\bf
V}_2=\mathrm{Im}(\phi_{\bar{g}_3})$.
 It is easy to check that  ${\bf V}_2$ is a
 two-dimensional  ${\mathbb F}_3-$algebra and  has no
 non-zero nilpotent elements. Moreover, there is a  direct product of
  ${\mathbb F}_5-$algebras
 $${\bf V}={\bf V}_1\oplus {\bf V}_2.$$

 Similarly, we can compute analogously as the following
 direct product of
 one-dimensional ${\mathbb F}_3-$algebras
$${\bf V}_2=\langle\bar{s}_2-\bar{s}_3\rangle\oplus
\langle\bar{s}_2+\bar{s}_3\rangle.$$ Thus, ${\bf
V}=\langle\bar{s}_1\rangle\oplus
\langle\bar{s}_2-\bar{s}_3\rangle\oplus
\langle\bar{s}_2+\bar{s}_3\rangle$. It is easy to check that
$$(2\bar{s}_1)^2=2\bar{s}_1,\;(\bar{s}_2-\bar{s}_3)^2=\bar{s}_2-\bar{s}_3,\;
(2(\bar{s}_2+\bar{s}_3))^2=2(\bar{s}_2+\bar{s}_3).$$ Set
$\bar{h}_1=2\bar{s}_1,\; \bar{h}_2=\bar{s}_2-\bar{s}_3,\;
\bar{h}_3=2(\bar{s}_2+\bar{s}_3).$ We have
$${\bf
V}=\langle\bar{h}_1\rangle\oplus \langle\bar{h}_2\rangle\oplus
\langle\bar{h}_3 \rangle.$$

Compute Gr\"obner basis $G^*_1$ of  $\{f,1-u\bar{h}_1\}$ in
${\mathbb F}_3[x,u]$ by
   Lemma 5.
   $G_1^*=[2+x^2+x, u+1],$ so $$ f_1=2+x^2+x.$$

 Similarly, we can compute that
 \begin{eqnarray*}
 G_2^*&=& [x^3+2x^2+1, 2+u],\\
 G_3^*&=& [x+1, t+1].
  \end{eqnarray*}
Therefore,
 \begin{eqnarray*}
 f_2&=& x^3+2x^2+1,\\
 f_3&=& x+1.
 \end{eqnarray*}
It yields the following factorization
$$f=f_1f_2f_3.$$

{\bf Remark:} The virtue of our approach to to compute primary
decomposition of zero-dimensional ideals over finite fields is
that it allows us to find all the primary components completely.
In particular, Proposition 1 contributes a new and simple method
to decompose the invariant subspace of the Frobenius map  on the
quotient algebra $k[{\bf x}]/I$ by theoretical and practical
considerations. But the complexity of our approach is mainly
depended  on computing Gr\"obner bases. There is no detailed
discussion of the complexity and implementation of our alternative
method here.



\begin{thebibliography}{99}




\bibitem{CM02}C. Monico, Computing the primary decomposition of
zero-dimensional ideals, {\it J. Symbolic Comput.} {\bf 34}
(2002), 451-459.

\bibitem{Berlek67}E.R. Berlekamp, ¡°Factoring polynomials over finite fields¡±,
Bell System Technical J., {\bf 46}(1967), 1853-1859.

\bibitem{NKY}N. Noro, K. Yokoyama, Implementation of Prime Decomposition of
Polynomial Ideals over Small Finite Fields, {\it J. Symbolic
Comp.} {\bf 38} (4), 2004, 1227-1246.




\bibitem{Cox:06}Cox, D., Little, J., O¡¯Shea, D., Ideal, Varieties,
and Algorithms. Springer-Verlag, 2006.

\bibitem{Cox98}Cox, D., Little, J., O¡¯Shea, D.,  Using algebraic geometry,
Springer-Verlag,  1998.

\bibitem{Martin:00}Martin K.,  and  Lorenzo, R.,  Computational Commutative Algebra 1. Springer-Verlag, 2000.



\bibitem{GWM2009}S. Gao, D. Wan and M. Wang, Primary Decomposition of Zero-dimensional Ideals over Finite Fields, {\it Mathematics of Computation}. {\bf 78} (2008), 509-521.





\bibitem{Shi96}T. Shimoyama and K. Yokoyama, Localization and primary
decomposition of polynomial ideals, {\it J. Symbolic Comput.} {\bf
22} (1996), 247-277.

\bibitem{GT88}Gianni, P., Trager, B., Zacharias, G., Gr¡§obner bases and
primary decomposition of polynomial ideals, {\it J. Symbolic
Comput.} {\bf 6} (1988),  149-167.











\end{thebibliography}
\end{document}